\documentclass[11pt, a4paper]{amsart}
\usepackage{color}
\usepackage{graphicx}
\usepackage[english,french]{babel}
\usepackage{amssymb}
\usepackage{amscd}
\usepackage{amsmath}
\usepackage{mathrsfs}
\def\RR{{\mathbb {R}}}

\def\ep{\varepsilon}

\def\di{\displaystyle}

\def\R{{\mathbb {R}}}
\def\uep{u_{\ep}}

\def\K{{\mathcal{K}}}
\def\J{{\mathcal{J}}}
\def\H{{\mathcal{H}}^{N-1}}

\def\pint{\operatorname {--\!\!\!\!\!\int\!\!\!\!\!--}}

\def\jep{{\mathcal J}_\ep}

\newtheorem{teo}{Theorem}[section]
\newtheorem{lema}{Lemma}[section]
\newtheorem{defi}{Definition}[section]

\newtheorem{corol}{Corollary}[section]
\theoremstyle{definition}
\newtheorem{remark}{Remark}[section]
\newtheorem*{ack}{Acknowledgements}
\renewcommand{\theequation}{\arabic{section}.\arabic{equation}}

\begin{document}

\selectlanguage{english}

\title[An optimization problem with volume constrain]
{An optimization problem with volume constrain in Orlicz spaces}

\author[S. Mart\'{i}nez]
{Sandra Mart\'{\i}nez}

\address{Departamento  de Matem\'atica, FCEyN
\hfill\break\indent UBA (1428) Buenos Aires, Argentina.}
\email{{\tt  smartin@dm.uba.ar } }

\thanks{Supported by ANPCyT PICT No.
03-13719 and 03-10608, UBA X052 and X066 and Fundaci\'on Antorchas
13900-5. }



\maketitle

\begin{abstract}
We consider the  optimization problem of minimizing
$\int_{\Omega}G(|\nabla u|)\, dx$ in the class of functions
$W^{1,G}(\Omega)$, with a constrain on the volume of $\{u>0\}$.
The conditions on the function $G$ allow for a different behavior
at
 $0$ and at $\infty$. We consider a penalization problem, and we prove that
for small values of the penalization parameter, the constrained
volume is attained. In this way we prove that every solution $u$
is locally Lipschitz continuous and that the free boundary,
$\partial\{u>0\}\cap \Omega$, is smooth.
\end{abstract}
\begin{otherlanguage}{french}
\begin{abstract}

 Nous consid\'{e}rons le probl\`{e}me d'optimisation de minimiser
 $\int_{\Omega}G(|\nabla u|)\, dx$ sur la classe des fonctions
 $W^{1,G}(\Omega)$, avec une restriction sur le volume de
 $\{u>0\}$. Les conditions sur la fonction $G$ permettent un
 comportement diff\'{e}rent en $0$ et \`{a} l\'infini. Nous
 consid\'erons un probl\`{e}me de p\'enalisation et nous prouvons
 que le volume fix\'{e} est atteint quand la valeur de la p\'{e}nalisation
 est petite. De cette mani\`{e}re nous prouvons que toute solution
 $u$ est localement Lipschitzienne et que la fronti\`{e}re libre
 $\partial\{u>0\}\cap \Omega$ est reguli\`{e}re.
 \end{abstract}
\end{otherlanguage}

\section{Introduction}

We  begin with a few historical remarks. In the paper \cite{AAC},
Aguilera, Alt and Caffarelli study an optimal design problem with
a volume constrain. The authors prove the  regularity of
minimizers by introducing a penalization term in the energy
functional (the Dirichlet integral) and minimizing without the
volume constrain.

The steps that they follow are the following. First, the authors
observe that, for fixed values of the penalization parameter, the
penalized functional is very similar to the one considered in the
paper \cite{AC}, then  the regularity results for minimizers of
the penalized problem follow almost without change as in
\cite{AC}. Finally, they prove that  for small values of the
penalization parameter, the constrained volume is attained. In
this way, all the regularity results apply to the solution of the
optimal design problem.

This method has been applied to other problems with similar
success. In \cite{ACS, FBRW, Le, Te}, where the differential
equation satisfied by the minimizers is nondegenerate, uniformly
elliptic and  in \cite{FBMW1},  where the equation involved  may
be degenerate or singular elliptic, but it steals has the property
of being  homogeneous.

In this article we show that the same kind of results can be
obtained if we study a problem such that the differential equation
satisfied by the minimizers is  nonlinear degenerate or singular
elliptic, and possibly not homogeneous. More precisely, the
operator here has the form $ \mathcal{
L}u=\mbox{div\,}\Big(g(|\nabla u|)\frac{\nabla u}{|\nabla
u|}\Big)$ where $g$ satisfies the natural conditions introduced by
Lieberman in \cite{Li1}. These conditions generalize the so called
natural conditions of Ladyzhenskaya and Ural'tseva. In \cite{Li1}
the author studies the regularity of weak solutions to the
equation
\begin{equation}\label{L}
{\mathcal L}u=0.
\end{equation}
Lieberman proves, that under these conditions,  solutions of
\eqref{L} are $C^{1,\beta}$.

The conditions imposed to
 $g$ are the following,
\begin{equation}\label{cond}
0<\delta\le\frac{t g'(t)}{g(t)}\le g_0\ \ \ \ \forall t>0
\end{equation}
for certain constants $\delta $ and $g_0$. Observe that
$\delta=g_0=p-1$ when $g(t)=t^{p-1}$, and conversely, if
$\delta=g_0$ then $G$ is a power. For more examples of functions
satisfying \eqref{cond} see \cite{MW1}.

Condition \eqref{cond} ensures that the equation \eqref{L} is
equivalent to a uniformly elliptic equation in nondivergence form
with ellipticity constants independent of the solution $u$ on sets
where $\nabla u\neq0$. This condition does not imply any kind of
homogeneity on the function $G$ (the primitive of $g$) and
moreover, it allows for a different behavior of the function $g$
when $|\nabla u|$ is close to zero or infinity.

We give  now, more precisely the description of the problem that
we study,

Take $\Omega$ a smooth bounded domain in $\R^N$ and $\varphi_0 \in
W^{1,G}(\Omega)$, a Dirichlet datum, with $\varphi_0\geq c_0>0$ in
$\bar{A}$, where A is a nonempty relatively open subset of
$\partial\Omega$ such that $A\cap\partial\Omega$ is $C^2$. Here
$W^{1,G}(\Omega)$ is a Sobolev-Orlicz space (see Appendix A). Let
$$
\K_\alpha = \{u\in W^{1,G}(\Omega)\, / \, |\{u>0\}| = \alpha, \,
u=\varphi_0\ \mbox{ on }\partial\Omega\}.
$$

Our problem is to minimize $\J(u)=\int_\Omega G(|\nabla u|)\,dx$
in $\K_\alpha$, with $g=G'$ satisfying \eqref{cond}.

One of the  difficulties of these problems is to prove the
regularity of the minimizers, since it is hard to make enough
volume preserving perturbations without the previous knowledge of
the regularity of $\partial\{u>0\}$.

In order to solve our original problem in a way that allows us to
perform non volume preserving perturbations we follow the idea of
\cite{AAC} and consider instead the following penalized problem:
We let
$$
\K=\{u\in W^{1,G}(\Omega)\, / \, u = \varphi_0\ \mbox{ on }
\partial\Omega\}
$$
and
\begin{equation}\label{jep}
\jep(u) = \int_\Omega G(|\nabla u|)\, dx + F_\ep(|\{u>0\}|),
\end{equation}
where
$$
F_\ep(s) = \begin{cases}
\ep(s-\alpha) & \mbox{ if } s<\alpha\\
\frac1\ep(s-\alpha) & \mbox{ if } s\ge\alpha.
\end{cases}
$$

Then, the penalized problem is
\begin{equation*}\label{Pep}
\tag{$P_\ep$} \mbox{Find } \uep\in \K \quad \mbox{such that} \quad
\J_\ep(\uep) = \inf_{v\in \K} \J_\ep(v).
\end{equation*}

In order to prove the  existence of minimizers we use some compact
immersion theorems in Sobolev-Orlicz spaces, and the result
follows easily by direct minimization. The regularity of the
minimizers and of their free boundaries $\partial \{ \uep
>0 \}$ follows by showing that any minimizer $\uep$ is a solution of the
following free boundary problem,
\begin{equation}\label{ecudebil}
\begin{cases}
\mathcal{L} u_{\ep}  = 0 &\mbox{ in } \{u_{\ep}>0\}\cap \Omega,\\
u_{\ep}=0,\quad \di\frac{\partial u_{\ep}}{\partial\nu} =
\lambda_\ep & \mbox{ on }
\partial\{u_{\ep}>0\}\cap\Omega,
\end{cases}
\end{equation}
in the sense that was defined in \cite{MW1}, where $\lambda_\ep$
is a positive constant. The properties of the definition of weak
solution are not difficult to establish since the minimization
problem studied in \cite{MW1} is very similar to \eqref{Pep}. The
only difference is that in \eqref{Pep} the functional is linear in
$|\{u>0\}|$ and  here the term $F_{\ep}$  is piecewise linear and
zero in $\alpha$. With these results we have that for almost $\H-$
every point, the free boundary is locally a $C^{1,\beta}$ surface
(see Corollary \ref{teo.regularity} in \cite{MW1}).

We also improve the regularity result for the case $N=2$, for a
subclass of functions satisfying \eqref{cond}. We prove, that in
this case, the whole free boundary is regular. Full regularity of
the free boundary in dimension 2 was prove  in \cite{AAC},
\cite{ACF} and in \cite{DP2} if $2-\delta\le p<\infty$ for a small
$\delta>0$. Also for the penalization problem in \cite{Le}.
 A similar result was
proved by A. Petrosyan in dimension 3 for $p$ close to 2 (see
\cite{P}).

As in \cite{AAC}, the reason why this penalization method is so
useful is that there is no need to pass to the limit in the
penalization parameter $\ep$ for which uniform, in $\ep$,
regularity estimates would be needed. In fact, we show that for
small values of $\ep$ the right volume is already attained. This
is, $|\{\uep>0\}|=\alpha$ for $\ep$ small. This step is where
 the proof parts from previous work on similar problems,
  since here we may not have the homogeneity
of the function $g$ (see Lemma \ref{lapapa}).

Finally, the fact that, for small $\ep$, any minimizer of $\J_\ep$
satisfies $|\{\uep>0\}|=\alpha$ implies that any minimizer of our
original optimization problem is also a minimizer of $\J_\ep$ so
that it is locally Lipschitz continuous with smooth free boundary.

\bigskip

The paper is organized as follows: In Section 2 we begin our
analysis of problem \eqref{Pep} for fixed $\ep$. First we prove
the existence of a minimizer, local Lipschitz regularity and
nondegeneracy near the free boundary (Theorem \ref{existencia})
and we prove that minimizers are weak solutions of  a free
boundary problem as defined in \cite{MW1} (Remark \ref{weak3}).
Then we have that for almost $\H-$ every   point, the free
boundary is locally a $C^{1,\beta}$ surface (Corollary
\ref{teo.regularity}). We prove that, for the case $N=2$, for a
subclass of functions satisfying \eqref{cond} their hole free
boundary is regular (Corollary \ref{regdim2}). In Section 3 we
prove that for small values of $\ep$ we recover our original
optimization problem.

We include at the end of the paper a couple of appendices where
some results about Orlicz spaces,   some properties of
$\mathcal{L}-$sub\-harmonic functions and blow up sequences are
established.


\section{The penalized problem}
\setcounter{equation}{0}

\subsection{Regularity of minimizers and their free boundaries }

We begin by discussing the existence of extremals and the
regularity. We are going to give some properties of the
minimizers. As the functional $\J_{\ep}$ is very similar to the
one in \cite{MW1}, some of the proof of these properties follows
as in \cite{MW1}. In that cases we are only  going to state the
results and avoid any proof. Next, we prove that any minimizer of
$\J_{\ep}$ is a weak solution of \eqref{ecudebil}, as was defined
in \cite{MW1}. Therefore we will have, by the results therein that
the free boundary is smooth.

\begin{teo}\label{existencia}
Let $\Omega\subset\R^N$ be bounded. Then there exists a solution
to the problem \eqref{Pep}. Moreover, any such solution $\uep$ has
the following properties:
\begin{enumerate}
\item $\uep$ is locally Lipschitz continuous in $\Omega$, and for
$D\subset\subset\Omega$ we have that, $\|\nabla
u\|_{L^{\infty}(D)}\leq C$ with $C=C(N,g_0,\delta,
dist(\partial\Omega,D),\ep)$. \item $\mathcal{L}u_{\ep}=0$ in
$\{u_{\ep}>0\}$.

\item There are constants $0< c_{min}\leq C_{max}$, $\gamma\ge1$,
such that for balls $B_r(x)\subset D$ with $x\in \partial
\{u_{\ep}>0\}$
$$
c_{min}\leq \frac{1}{r}\Big(\pint_{B_r(x)} u_{\ep}^{\gamma} dx\,
\Big)^{1/\gamma}\leq C_{max}
$$

\item For every $D\subset\subset\Omega$, there exist constants $C,
c>0$ such that for every $x\in D\cap \{\uep>0\}$,
$$
c\, {\rm dist}(x,\partial\{\uep>0\})\le \uep(x)\le C\, {\rm
dist}(x,\partial\{\uep>0\}).
$$

\item For every $D\subset\subset\Omega$, there exists a constant
$c>0$ such that for $x\in \partial\{u_{\ep}>0\}$ and
$B_r(x)\subset D$,
$$
c\le \frac{|B_r(x)\cap\{\uep>0\}|}{|B_r(x)|}\le 1-c.
$$
\end{enumerate}
The constants may depend on $\ep$.
\end{teo}

\begin{proof}
Observe that, if $A\leq B$ then $\ep (B-A)\leq F_{\ep}(B)-
F_{\ep}(A)\leq \frac{1}{\ep} (B-A)$. Then the proof follows as in
section 3, 4 and 5 in \cite{MW1}.

\end{proof}

From now on we drop the subscript $\ep$ and denote by $u$ instead
of $\uep$ a solution to \eqref{Pep}.

\begin{teo}[Representation Theorem] \label{DP-teo5.2} Let $u\in \K$ be a
solution to \eqref{Pep}. Then,
\begin{enumerate}
\item $\H( D\cap\partial\{u>0\})<\infty$ for every
$D\subset\subset\Omega$.

\item There exists a Borel function $q_u$ such that
$$
\mathcal{L}u= q_u \,\H \lfloor \partial\{u>0\}.
$$

\item For $D\subset\subset\Omega$ there are constants $0<c\le
C<\infty$ depending on $N, \Omega, D$ and $\ep$ such that for
$B_r(x)\subset D$ and $x\in \partial\{u>0\}$,
$$
c\le q_u(x)\le C,\quad c\,r^{N-1}\le
\H(B_r(x)\cap\partial\{u>0\})\le C\,r^{N-1}.
$$

\item $\H(\partial\{u>0\}\setminus\partial_{\rm red}\{u>0\})=0.$
\end{enumerate}
\end{teo}
\begin{proof}
For the proof, see sections 6 and 7 in \cite{MW1}. Observe that
$D\cap\partial\{u>0\}$ has finite perimeter, thus, the reduce
boundary $\partial_{\rm red}\{u>0\}$ is defined as well as the
measure theoretic normal $\nu(x)$ for $x\in\partial_{\rm
red}\{u>0\}$ (see \cite{F}).
\end{proof}

\begin{lema}\label{blow} Let $x_0, x_1\in \partial\{u>0\}$ and
$\rho_k\to 0^+$. For $i=0,1$ let $x_{i,k}\to x_i$ with
$u(x_{i,k})=0$ such that $B_{\rho_k}(x_{i,k}) \subset \Omega$ and
such that the blow-up sequence
$$
u_{i,k}(x) = \frac{1}{\rho_k} u(x_{i,k} + \rho_k x)
$$
has a limit $u_i(x)=\lambda_i ( x\cdot \nu_i)^-$, with
$0<\lambda_i<\infty$ and $\nu_i$ a unit vector. Then $\lambda_0 =
\lambda_1$.
\end{lema}
\begin{proof}
It follows as in \cite{FBMW1} by using the results in Appendix C.
\end{proof}

\begin{lema}\label{limsup1}Let $x_0\in\Omega\cap\partial\{u>0\}$ and let
$$
\lambda=\lambda(x_0):=\limsup_{\stackrel{x\to x_0}{u(x)>0}}
|\nabla u(x)|.
$$
Then there exists  sequences $y_k$, $d_k$ and $\nu_k,\nu$ such
that $\nu_k\rightarrow \nu$ and the blow up sequence with respect
to $B_{d_k}(y_k)$ has limit,
$$
u_0(x)=\lambda   (x\cdot \nu)^-.
$$
\end{lema}
\begin{proof}
It follows as in the proof of Theorem 2.3 in \cite{FBMW1} by using
the results in Appendix B and C.
\end{proof}

\begin{lema}\label{q1}
For $\H$--a.e. $x_0\in \partial_{\rm red}\{u>0\}$, there exist a
sequence $\gamma_n$  such that if $u_n$ is  the blow up sequence
with respect to $B_{\gamma_n}(x_0)$ we have that,
$$
u_n\rightarrow {\lambda^*}   (x\cdot \nu(x_0))^-
$$
with $\nu(x_0)$ the outward unit normal to $\partial\{u>0\}$ in
the measure theoretic sense and ${\lambda^*}=g^{-1}(q_u(x_0)).$
\end{lema}
\begin{proof}
Suppose that $\nu(x_0)=e_N$. As in Theorem 3.5 in \cite{ACF} and
Theorem 5.5 in \cite{DP} we can prove by using the boundary
regularity of solutions of $\mathcal{L}v=0$ (see \cite{Li1}) that
for $\H$--a.e. $x_0\in
\partial_{\rm red}\{u>0\}$,  the following fact holds. If we consider
the blow up limit $u_0$ of $u$ with respect to sequences of balls
$B_{\rho_k}(x_0)$, $\rho_k\rightarrow 0$ we have that,
\begin{equation}\label{cauchy}
\begin{cases}
\mathcal{L}u_0=0 &\quad \mbox{in} \{x_N<0\}\\
u_0=0, \quad g(|\nabla u_0|)=q_u(x_0) &\quad \mbox{on} \{x_N=0\}.
\end{cases}
\end{equation}
Therefore, $u_0(x)={\lambda^*}  x_N^-+o(|x|)$.

 Take now $u_{0,j}$, a blow up sequence of $u_0$,
 with respect to balls $B_{\mu_j}(0)$, therefore
$$u_{0,j}\rightarrow u_{00}={\lambda^*}  x_N^-.$$
Now, we want to construct a blow up sequence of $u_0$ with limit
$u_{00}$. Observe, that
$$\Big|\frac{1}{\rho_k \mu_j}u(x_0+\rho_k \mu_j x)-u_{00}(x)\Big|\leq
\frac{1}{\mu_j}|u_k(\mu_jx)-u_{0}(\mu_j
x)|+|u_{0,j}(x)-u_{00}(x)|,
$$ and since
$u_{k}\rightarrow u_0$ uniformly over compacts sets we have that
for $j\geq j_n$, $|u_{0,j}(x)-u_{00}(x)|<1/n$ and for $k\geq
k_{j,n}$, $|u_k(\mu_jx)-u_{0}(\mu_j x)|<\mu_j/n $ if $|x|\leq n$.
We may suppose that $j_n\geq n$ and $k_{j,n}\geq n$. Taking
$j=j_n$, $k=k_{j_n,n}$, and $\gamma_n=\rho_{k_{j_n,n}} \mu_{j_n}$.
We have that $\gamma_n\to 0$ and $|u_{\gamma_n}(x)-u_{00}(x)|<2/n$
in $B_n$. The result follows.
\end{proof}
\begin{teo}\label{const}
Let $u\in \K$ be a solution to \eqref{Pep} and $q_{u}$ the
function in Theorem \ref{DP-teo5.2}. Then there exists a constant
$\lambda_u$ such that
\begin{align}
& \limsup_{\stackrel{x\to x_0}{u(x)>0}} |\nabla u(x)| =
\lambda_u,\qquad \mbox{for every } x_0\in
\Omega\cap\partial\{u>0\}\label{limsup}\\
& q_{u}(x_0)=g(\lambda_u),\qquad \H -\mbox{a.e } x_0\in
\Omega\cap\partial_{red}\{u>0\}. \label{qu}
\end{align}
\end{teo}
\begin{proof}
It follows as in \cite{Le} by using Lemmas \ref{blow},
\ref{limsup1}, \ref{q1} and Theorem \ref{DP-teo5.2} (4).
\end{proof}
Now, we can prove the asymptotic development for minimizers,
\begin{teo}\label{asimp}
For $\H-$a.e. point $x_0\in\partial\{u>0\}$ belongs to
$\partial_{red}\{u>0\}$ and
$$
u(x_0+x) = \lambda_{u} (x\cdot \nu(x_0))^- + o(|x|) \quad
\mbox{for}\quad x\to 0.
$$
\end{teo}
\begin{proof}
We can prossed as in the proof of Lemma \ref{q1} until arrive to
equation \eqref{cauchy}. Now, as by Theorem \ref{const} we have
that $q_u(x_0)=g(\lambda_u)$ and $|\nabla u|\leq \lambda_u$, we
can use the same argument of Theorem 5.5 in \cite{DP} and conclude
that $u_0(x)={\lambda_u}  x_N^-$. And as, the blow up sequence was
arbitrary chosen, we have the desired result.
\end{proof}
\medskip

\begin{remark}\label{weak3}Now we have, by properties (1), (2), (3) in Theorem
\ref{existencia} and Theorem \ref{asimp} that any minimizer
satisfies all the properties of the definition of weak solution II
in \cite{MW1}. Therefore we have by Theorem 9.3 and Remark 9.2 in
\cite{MW1} and Theorem \ref{DP-teo5.2} (4) the following
regularity result for the free boundary $\partial \{u>0\}$.
\end{remark}
\begin{corol}\label{teo.regularity}
Let $u\in \K$ be a solution to \eqref{Pep}. Then there exists a
subset $A\subset\partial_{red}\{u>0\}$ with $\H
(\partial_{red}\{u>0\}\setminus A)=0$ such that for any $x_0\in A$
there exists $r>0$ so that $B_r(x_0)\cap\partial\{u>0\}$ is a
$C^{1,\alpha}$ surface. The remainder of the free boundary has
$\H-$measure zero.
\end{corol}

\subsection{ Full regularity for the case $N=2$ }
We will prove, that in dimension two, for a subclass of functions
satisfying \eqref{cond},  their hole free boundary is a
$C^{1,\beta}$ surface.

The class that we consider consists on those functions satisfying
condition \eqref{cond} and such that,
\begin{equation}\label{condi}
\mbox{There exist  constants }  t_0>0 \mbox{ and } k>0 \mbox{ so
that } g(t)\leq k t \mbox{ for } t\leq t_0.
\end{equation}

Observe that this condition is satisfied for example, if
$\delta\geq 1$ or when
 $g_0\geq 1$ and there exists a constant $C$ such that
$\displaystyle\lim_{t\to 0} \frac{g(t)}{t^{g_0}}=C.$

In order to prove the full regularity, we first need the following
two Lemmas, that hold for any dimension and for any $\delta$ and
$g_0$,

\begin{lema}\label{estimacion.del.gradiente}
Let $u\in \K$ be a local minimizer. Given $D\subset\subset\Omega$,
there exist constants $C=C(N,D,\lambda_u)$, $r_0=r_0(N,D)>0$ and
$\gamma=\gamma(N,D)>0$ such that, if $x_0\in D\cap\partial\{u>0\}$
and $r<r_0$, then
$$\sup_{B_r(x_0)} |\nabla u|\le \lambda_u+Cr^\gamma.$$
\end{lema}

\begin{proof} The proof is similar to the proof of Theorem 7.1 in
\cite{DP} but here we make a little modification by using a result
of \cite{Li2}. This result allows us not having to add any new
hypothesis to the function $g$.

Let $U_{\rho}=\big(G(|\nabla u|) - G(\lambda_u) - \rho\big)^+$ and
$U_0=\big(G(|\nabla u|) - G(\lambda_u)\big)^+$. By Theorem
\ref{const} we know that $U_{\rho}$ vanishes in a neighborhood of
the free boundary. Since $U_{\rho}>0$ implies $G(|\nabla
u|)>G(\lambda_u)+\rho$, the closure of $\{U_{\rho}>0\}$ is
contained in $\{G(|\nabla u|)>G(\lambda_u)+\rho/2\}$.
 The function $u$ satisfies the linearized equation
 $$Tu=b_{ij}(\nabla u) D_{i j}u=0$$ where
 $b_{ij}$ is defined in \eqref{opT}, and is
 $\beta$-elliptic in $\{G(|\nabla
u|)>G(\lambda_u)+\rho/2\}$.

Let $v=G(|\nabla u|)$, by Lemma 1 in \cite{Li2}, we have that $v$
satisfies,
$$Mv=D_i(b_{ij}(\nabla u) D_j v)\geq 0 \ \ \mbox{ in
}\{G(|\nabla u|)>G(\lambda_u)+\rho/2\}.$$ Hence $U_{\rho}$
satisfies
$$MU_{\rho}\geq 0 \ \ \mbox{ in } \{G(|\nabla u|)>G(\lambda_u)+\rho/2\}.$$
Extending the operator $M$ with the uniformly elliptic
divergence-form operator
$$\widetilde{M}w=D_i(\widetilde{b}_{ij}(x) D_j w)\ \ \mbox{ in }
\Omega$$ with measurable coefficients such that
$$\widetilde{b}_{ij}(x)= b_{ij}(\nabla u)
 \ \ \mbox{ in } \{G(|\nabla u|)>G(\lambda_u)+\rho/2\},$$
we obtain
$$\widetilde{M}U_{\rho}\geq 0\ \ \mbox{ in } \Omega.$$
For any $r>0$ set
$$
h_{\rho}(r) = \sup_{B_r(x_0)} U_{\rho},\qquad h_0(r) =
\sup_{B_r(x_0)} U_0,
$$
for any $r< r_0=\,{\rm dist}\,(D,\partial\Omega)$ and $x_0\in
D\cap\partial\{u>0\}$.

Then, $h_{\rho}(r) - U_{\rho}$ is a $\widetilde{M}$- supersolution
in the ball $B_r(x_0)$ and
$$
\begin{array}{lll}
h_{\rho}(r) - U_{\rho} & \ge 0 & \mbox{in } B_r(x_0)\\
 & = h_{\rho}(r) & \mbox{in } B_r(x_0)\cap\{u=0\}.
\end{array}
$$
Applying the weak Harnack inequality (see \cite{GT} Theorem 8.18)
with $1\le p< N/(N-2)$, we get
$$
\inf_{B_{r/2}(x_0)} \big(h_{\rho}(r) - U_{\rho}\big)\ge c r^{-N/p}
\|h_{\rho}(r) - U_{\rho}\|_{L^p(B_r(x_0))}\ge c h_{\rho}(r),
$$
since, by Theorem \ref{existencia}, $|B_r(x_0)\cap\{u=0\}|\ge c
r^N$. Taking now $\rho\to 0$ we obtain
$$
\inf_{B_{r/2}(x_0)} \big(h_0(r) - U_0\big)\ge c h_0(r),
$$
for some $0<c<1$, which is the same as
$$\sup_{B_{r/2}(x_0)} U_0 \le (1-c) h_0(r).$$
Therefore
$$h_0\left(\frac{r}{2}\right)\le (1-c)h_0(r),$$
from which it follows that $h_0(r)\le C r^{\gamma}$ for some
$C>0$, $0<\gamma<1$. That is,
$$G(|\nabla u|)\leq G(\lambda_u)+ C r^{\gamma}$$ and
therefore
$$|\nabla u|\leq \lambda_u+ C r^{\gamma}$$ and
 now the conclusion of the Theorem follows.
\end{proof}

\begin{lema}\label{promphi1}
Let $x_1$ be a regular   free boundary point.

Take
$$
\tau_{\rho}(x)=\begin{cases} \displaystyle
x+\rho^2\phi\left(\frac{|x-x_1|}{\rho}\right)
\nu_{u}(x_1) & \mbox{ for } x\in B_{\rho}(x_1),\\
x &\mbox{ elsewhere, }
\end{cases}
$$
where $\phi\in C_0^{\infty}(-1,1)$ with $\phi'(0)=0$.

Let
\begin{equation}\label{deltarho}
\delta=\rho^2 \int_{B_{\rho}(x_1)\cap\partial\{u>0\}}
\phi\left(\frac{|x-x_1|}{\rho}\right)\, d\H.
\end{equation}
Take $v_{\delta}(x)=v_\rho(x)=u(\tau_{\rho}^{-1}(x))$, then
\begin{align}\label{chica2}
\int_{B_{\rho}(x_1)}(G(|\nabla v_\rho|) - G(|\nabla u|))\, dx = -l
\rho^{N+1} \Phi(\lambda_u) + o(\rho^{N+1}),
\end{align}
where $l=\lim_{\rho\to 0}\frac{\delta}{\rho^{N+1}}$ and
$\Phi(t)=g(t)t-G(t)$.

\end{lema}
\begin{proof}
The proof follows the lines of Theorem 3.1 in \cite{FBMW1}.
\end{proof}

In the following Lemma is where we need to impose condition
\eqref{condi}.

\begin{lema}\label{promphi2}
Let $\Phi(t)=g(t) t-G(t)$, and $g$ satisfying condition
\eqref{condi}. Let $x_0$ be a free boundary point ,
$D\subset\subset \Omega$ and $B_{\mu}(x_0)\subset D$. Take
$v=\max(u-t\eta,0)$, where $t>0$, $\eta\in C_0^{\infty}(\Omega)$,
$\eta=0$ in $\Omega\setminus B_{\mu(x_0)}$ and $|\nabla \eta|\leq
C/t$. Therefore,
$$
\int_{B_{\mu}(x_0)\cap\{u>0\}}(G(|\nabla v|) - G(|\nabla u|))\,
dx\leq \int_{B_{\mu}(x_0)\cap\{0<u\leq t\eta\}} \Phi(|\nabla u|)\,
dx+ C_0 t^2 \int_{B_{\mu}(x_0)\cap\{u> t\eta\}}|\nabla \eta|^2\,
dx,
$$ for $C_0=C_0(N,\delta,g_0,dist(\partial \Omega,D),\ep, C)$.
\end{lema}
\begin{proof}
The Lemma follows as in Theorem 4.3 in \cite{ACF}. We only have to
make the following observations. First, observe that $|\nabla
u-t\nabla \eta|\leq |\nabla u|+C\leq C_1+C$, where $C_1$ is the
constant in Theorem \ref{existencia} (1). On the other hand, if
$g$ satisfies \eqref{condi}, and if $F(s)=\frac{g(s)}{s}$ then for
$0\leq s\leq C_1+C$, there exists a constant $C_0$ such that
$F(s)\leq C_0$.  Therefore we have that $F(|\nabla u-t\nabla
\eta|)$ is bounded by $C_0$. The rest of the proof follows as in
\cite{ACF}.
\end{proof}

Now, following ideas of \cite{Le}, using Lemmas
\ref{estimacion.del.gradiente},  \ref{promphi1} and
\ref{promphi2}, we prove, for $N=2$ and $g$ satisfying
\eqref{condi} the following,
\begin{teo}\label{promediophi}
Let $N=2$,  $g$ satisfying \eqref{condi} and $u$  a minimizer,
then for any ball $B_r$ centered at the free boundary we have,
$$\pint_{B_r\cap\{u>0\}} (\Phi(\lambda_u)-\Phi(|\nabla
u|))^+\rightarrow 0 \mbox{ when } r\to 0,$$ where $\Phi(t)=g(t)
t-G(t)$.
\end{teo}
\begin{proof}
Let $0<r<\mu$,  $t>0$ and $v_0$ be the function defined in Lemma
\ref{promphi2}. By Theorem \ref{existencia} $u\leq Cr$ in
$B_{r}(x_0)$, take $t=Cr$ and let $\delta_t=|\{0<u\leq t\eta\}\cap
B_{\mu}(x_0)|$.

Now, let us take $x_1$  far from $x_0$ and such that
$\partial\{u>0\}\cap B_{r_1}(x_1)$ is regular, for $r_1$ small.
Let $\rho$ be such that \eqref{deltarho} is satisfied for
$\delta=\delta_t$, and consider $v_1=v_{\delta_t}$ defined in
$B_{r_1}(x_1)$ as in Lemma \ref{promphi1}. Then, the function
$$v=\begin{cases}v_0 \quad\mbox{ in } B_{\mu}(x_0)
 \\ v_1 \quad \mbox{ in } B_{r_1}(x_1) \\ u \quad\mbox{ elsewhere } \end{cases}$$
is admissible for our minimization problem and
$|\{v>0\}|=|\{u>0\}|$. Therefore, by Lemmas \ref{promphi1} and
\ref{promphi2} we have
$$\begin{aligned} 0&\leq
\mathcal{J}_{\ep}(v)-\mathcal{J}_{\ep}(u)=\int_{B_{\rho}(x_0)}
(G(|\nabla v|)- G(|\nabla u|))\, dx +\int_{B_{r_1}(x_1)}
(G(|\nabla v|)- G(|\nabla u|))\, dx\\ &\leq \int_{
B_{\mu}(x_0)\cap\{u\leq t \eta\}} \Phi(|\nabla u|) + Ct^2
\int_{B_{\mu}(x_0)\cap\{u>t\eta\}}|\nabla \eta|^2\, dx -l \rho^{3}
\Phi(\lambda_u)+o(\rho^{3}).\end{aligned}$$ By definition of
$\delta_t$ we have,
$$\int_{B_{\mu}(x_0)\cap\{0<u\leq t \eta\}}
(\Phi(\lambda_u)-\Phi(|\nabla u|))\, dx \leq  Ct^2
\int_{B_{\mu}(x_0)\cap\{u>t\eta\}}|\nabla \eta|^2\, dx
+o(\rho^{3})+(\delta_t-l \rho^{3}) \Phi(\lambda_u).$$
 Now choose
$$\eta(x)=\begin{cases} \frac{\log(\mu/|x-x_0|)}{\log(\mu/r)}
&\quad \mbox{ in } B_{\mu}(x_0)\setminus B_{r}(x_0),\\1 &\quad
\mbox{ in } B_{r}(x_0)\\
0 &\quad \mbox{ in } \Omega\setminus B_{\mu}(x_0),\end{cases}$$
observe that the condition $|\nabla \eta|\leq C/t$ is satisfied if
we choose $\mu\geq 2 r$.

 By our election of $t$ and $\eta$
we have,
$$\begin{aligned}\int_{B_{r}(x_0)\cap\{u>0\}}
(\Phi(\lambda_u)-\Phi(|\nabla u|))^+\, dx &\leq
\int_{B_{\mu}(x_0)} (\Phi(|\nabla
u|)-\Phi(\lambda_u))^+\, dx + \frac{Cr^2}{\log(\mu/r)}\\
&+o(\rho^{3})+(\delta_t-l \rho^{3})
\Phi(\lambda_u).\end{aligned}$$ By Lemma
\ref{estimacion.del.gradiente}, we have that $\Phi(|\nabla
u|)-\Phi(\lambda_u)\leq \Phi(\lambda_u+C
r^{\gamma})-\Phi(\lambda_u)=\Phi'(\xi) Cr^{\gamma}$, for
$\lambda_u\leq \xi \leq \lambda_u+C r^{\gamma}$. As
$\Phi'(t)=g'(t)t\leq g_0 g(t)$, and as $g$ is increasing we have
that $\Phi'(\xi)\leq g_0 g(\xi)\leq g_0 g(\lambda_u+C
r^{\gamma})$.

Therefore by definition of $l$ we have
$$\begin{aligned}\pint_{B_{r}(x_0)\cap\{u>0\}}
(\Phi(\lambda_u)-\Phi(|\nabla u|))^+\, dx &\leq C\Big(
\frac{(\mu^{\gamma+2}+o(\rho^3))}{r^2}+
\frac{1}{\log(\mu/r)}\Big),\end{aligned}$$ where $C=C(\lambda_u)$.
As by Theorem \ref{existencia} (5), $\delta_t\leq c \mu^2$ we have
that $o(\rho^3)=o(\mu^2)$. Taking $r=\mu h(\mu)^\beta$, where
$h(\mu)=\max\Big(\mu, \frac{o(\mu^2)}{\mu^2}\Big)$ with
$\beta<\min\{\gamma/2, 1/2\}$, we have the desired result.
\end{proof}
\begin{corol}\label{regdim2}
Let $N=2$, $g$ satisfying \eqref{condi} and $u\in \K$ be a
solution to \eqref{Pep}. Then $\partial\{u>0\}$ is a $C^{1,\beta}$
surface locally in $\Omega$.
\end{corol}
\begin{proof}
The proof follows now as in \cite{AC}, we give the proof here for
the readers convenience. Let $u_k$ be a blow up sequence
converging to $u_0$. Since,  $\nabla u_k\rightarrow \nabla u_0$
a.e in $\mathbb{R}^N$, we conclude from Theorem \ref{const} and
Theorem \ref{promediophi} that $|\nabla u_0|=\lambda_u$ in
$B_1\cap\{u_0>0\}$, and then
$$0=\mathcal{L}u_0=div\Big(\frac{g(|\nabla u_0|)}{|\nabla u_0|}\nabla u_0\Big)=
\frac{g(\lambda_\ep)}{\lambda_u} \triangle u_{0}\quad \mbox{ in }
\{u_0>0\}.$$ Therefore $u_0$ is harmonic in $\{u_0>0\}$, and if we
take $v=|\nabla u_0|^2$, we have $0=\triangle v=|D^2u_0|^2$ and
that means that $\nabla u_0$ is constant in each connected
component of this set.  Therefore, by Lemma C.1 (6) and (8) we
have,
$$u_0= \lambda_u \max(x\cdot \nu_0,0) +q \max(-x\cdot
\nu_0,s),$$ for some $\nu_0$ and  $q, s\geq 0$. Since $\{u_0=0\}$
has positive density at the origin, we have that $s>0$ or $q=0$.
Therefore, we have proved that any blow up sequence has a
subsequences that converges to a half linear function $u_0=
\lambda_u \max(x\cdot \nu_0,0)$ in some neighborhood of the
origin, then applying Theorem 9.3 and Remark 9.2 in \cite{MW1} we
have the desired result.
\end{proof}

\begin{remark} Since the functional in \cite{MW1} is linear in
$|\{u>0\}|$ we can also prove, for minimizers of that problem, the
full regularity of the free boundary when $N=2$ . We only have to
use Theorem \ref{estimacion.del.gradiente}, Lemma \ref{promphi2}
(to treat the first term of the functional) and finally the result
follows as in \cite{AC}.
\end{remark}

\section{Behavior of the minimizer for small $\ep$.}
\setcounter{equation}{0} \label{sect.ep} \setcounter{equation}{0}
In this section, since we want to analyze the dependence of the
problem with respect to $\ep$ we will again denote by $\uep$ a
solution to problem \eqref{Pep}.

To complete the analysis of the problem, we will now show that if
$\ep$ is small enough, then
$$
|\{\uep>0\}| = \alpha.
$$
To this end, we need to prove that the constant $\lambda_\ep :=
\lambda_{\uep}$ is bounded from above and below by positive
constants independent of $\ep$. We perform this task in a series
of lemmas.

\begin{lema}\label{AAC-lema5}
Let $\uep\in \K$ be a solution to \eqref{Pep}. Then, there exists
a constant $C>0$ independent of $\ep$ such that
$$
\lambda_\ep  \le C.
$$
\end{lema}

\begin{proof}
The proof is similar to the one in Theorem 3 in \cite{AAC}.

First we will prove that there exist $C,c>0$, independent of
$\ep$, such that
$$
c\leq |\{\uep>0\}|\leq C\ep+\alpha.
$$
Taking $u_0$ such that $|\{u_0>0\}|\leq \alpha$ we have that
$\J_{\ep}(u_0)\leq C$ then we have that $F_{\ep}(|\{\uep>0\}|)
\leq C$ thus obtaining the bound from above. We also have that
$\int_{\Omega}G(|\nabla u_{\ep}|)$ is bounded. As  $u_{\ep} =
\varphi_0$ in $\partial\Omega$, we have by Lemma \ref{equi} that
$\|\nabla u_{\ep}-\nabla \varphi_0\|_{G}\leq C$ and by
Lemma~\ref{poinc} we also have $\|u_{\ep}-\varphi_0\|_{G}\leq C$,
then $\|\uep\|_{W^{1,G}(\Omega)} \leq C$. Using the Sobolev trace
Theorem, the H\"{o}lder inequality and the embedding Theorem
\ref{imb} we have, for $q<\delta+1$
$$
\begin{aligned}\int_{\partial \Omega} \varphi_0^{q}\, d\H &\leq C
|\{\uep>0\}|^{\frac{\delta+1-q}{\delta+1}}\|\uep\|_{W^{1,\delta+1}(\Omega)}^{q}\leq
C |\{\uep>0\}|^{\frac{\delta+1-q}{\delta+1}}
\|\uep\|_{W^{1,G}(\Omega)}^{q}\\ &\leq C
|\{\uep>0\}|^{\frac{\delta+1-q}{\delta+1}} , \end{aligned}$$ and
thus we obtain the bound from below.

The rest of the proof follows as in Lemma 3.1 in \cite{FBMW1}.

\end{proof}

\begin{lema}\label{gamapromedio}
Let $\uep\in \K$ be a solution to \eqref{Pep}, $B_r \subset\subset
\Omega$  and $v$ a solution to
$$
\mathcal{L} v=0\quad\mbox{in }B_r,\quad \quad v=\uep\quad\mbox{on
}\partial B_r.
$$
Then there exists a positive constant
$\gamma=\gamma(\delta,g_0,N)$ such that
$$
\int_{B_r} |\nabla (u_{\ep}-v)|^q\, dx \geq C |B_r \cap
\{u_{\ep}=0\}| \left(\frac{1}{r}\left(\di\pint_{B_r}
u_{\ep}^{\gamma}\, dx \right)^{1/\gamma}\right)^q
$$
for all $q\geq 1$ and where $C$ is a constant independent of
$\ep$.

\end{lema}

\begin{proof}
The proof follows the lines of Lemma 3.2 in \cite{FBMW1}. The only
difference here is that we have to use the weak Harnack inequality
 of \cite{Li1} (Theorem 1.3) instead of the known one.

\end{proof}

Without losing generality, from now on we will suppose that
$g_0\geq 1$.

\begin{lema}\label{lapapa}
Let $\uep$ and $v$ be as in  Lemma \ref{gamapromedio}, then if $r$
is small enough (depending on $\ep$) we have,
\begin{equation}\begin{aligned}\label{rn5}
\int_{B_r} (G(|\nabla u|)-G(|\nabla v|))\, dx \geq C\int_{B_r}
|\nabla \uep-\nabla v|^{g_0+1}\, dx
\end{aligned}\end{equation}
for some constant $C$ independent of $\ep$.

\end{lema}
\begin{proof}
First we will use an inequality  proved in \cite{MW1} (see Theorem
2.3). Let,
$$A_1=\{x\in {B_r}: |\nabla \uep-\nabla v|\leq 2 |\nabla \uep|\},\quad
A_2=\{x\in {B_r}: |\nabla \uep-\nabla v| > 2 |\nabla \uep|\},$$
then ${B_r}=A_1\cup A_2$ and we have that,
\begin{equation}\label{a1a2}
\int_{B_r} (G(|\nabla \uep|)-G(|\nabla v|))\, dx \geq C\Big(\int_{
A_2} G(|\nabla \uep-\nabla v|)  \, dx+\int_{A_1}F(|\nabla \uep|)
|\nabla \uep-\nabla v|^2 \, dx\Big).
\end{equation}
Therefore we have, using that $g_0\geq 1$ and  (g1) in Lemma
\ref{prop}, that when $|\nabla \uep|\leq 1$ and $|\nabla v-\nabla
u_{\ep}|\leq 1$,
\begin{equation}\label{vamos}\begin{aligned} G(|\nabla u_{\ep}-\nabla v|)&\geq C |\nabla \uep-\nabla
v|^{g_0+1}\\
 F(|\nabla \uep|)&\geq C |\nabla \uep|^{g_0-1}\geq C |\nabla \uep-\nabla
 v|^{g_0-1}\quad \mbox{ in } A_1.
\end{aligned}\end{equation}
On the other hand, by Lemma \ref{AAC-lema5} and \eqref{limsup}, we
have that for small $r$ (depending on $\ep$), $|\nabla \uep|$ is
bounded by a constant independent of $\ep$. By Lemma 5.1 in
\cite{Li1} we have that there exist
$C_0,C_1=C_0,C_1(N,g_0,\delta)$ such that,
$$\sup_{B_r} G(|\nabla v|)\leq \frac{C_0}{r^N} \int_{B_{2r}}
G(|\nabla v|)\, dx \leq \frac{C_1}{r^N} \int_{B_{2r}} (1+G(|\nabla
u_{\ep}|))\, dx \leq \bar{C},
$$
if we choose $r$ small (depending on $\ep$) and where $\bar{C}$ is
independent of $\ep$. Then $|\nabla u_{\ep}|$ and $|\nabla
u_{\ep}-\nabla v|$ are bounded in $B_r$ by a constant independent
of $\ep$.
 Therefore,
\eqref{vamos} holds for all $x\in B_r$ and for a constant $C$
(independent of $\ep$). Combining \eqref{a1a2} and \eqref{vamos}
we obtain the desired result.

\end{proof}

\begin{lema}\label{pmayor2}
For every $\ep>0$ there exists a neighborhood of $A$ in $\Omega$
such that $\uep>0$ in this neighborhood.
\end{lema}
\begin{proof}
The proof follows the lines of the one in Lemma 3.4 in
\cite{FBMW1}. There is one step that it is convenient to mention
here. When we use the Schwartz symmetrization, we have to use that
this symmetrization preserves the distribution function and
strictly decreases the functional $\int_B G(|\nabla u|)\, dx$,
unless the function is already radially symmetric and radially
decreasing. These facts holds by Corollary 2.35,  in section II.8
of \cite{K}. The rest of the proof follows without any change.
\end{proof}
\begin{lema}\label{cotaabajo}
Let $\uep\in\mathcal{K}$ be a solution to $(P_{\ep})$, then
$$
\lambda_{\ep}\geq c>0,
$$
where $c$ is independent of $\ep$
\end{lema}

\begin{proof}
The proof follows as in \cite{FBMW1} by using Lemmas
\ref{gamapromedio}, \ref{lapapa},  \ref{pmayor2} and Lemma
\ref{propblowup} .
\end{proof}

With these uniform bounds on $\lambda_{\ep}$, we can prove the
desired result.
\begin{teo}\label{final}
Under the same hypotheses of Lemma \ref{cotaabajo}, there exists
$\ep_0>0$ such that for $\ep<\ep_0$, $|\{u_{\ep}>0\}|=\alpha$.
Therefore, $u_{\ep}$ is a minimizer of $\J$ in $\K_{\alpha}$.
\end{teo}
\begin{proof}
It follows as in Theorem 3.1 in \cite{FBMW1} be using Lemmas
\ref{AAC-lema5} and  \ref{cotaabajo}.
\end{proof}

As a corollary, we have the desired result for our problem

\begin{corol}
Any minimizer  $u$ of $\J$ in $\K_{\alpha}$ is a locally Lipschitz
continuous function,  $\partial_{red}\{u>0\}$ is a $C^{1,\beta}$
surface locally in $\Omega$ and the remainder of the free boundary
has $\H-$measure zero. Moreover if  $N=2$ and $g$ satisfies
\eqref{condi} then $\partial\{u>0\}$ is a $C^{1,\beta}$ surface
locally in $\Omega$.
\end{corol}
\begin{proof}
If $u$ is minimizer of $\J$ in $\K_{\alpha}$, by  Theorem
\ref{final} we have that for small $\ep$ there exists a solution
$\uep$ to $(P_{\ep})$ such that $|\{\uep>0\}|=\alpha$, then $u$ is
a solution to $(P_{\ep})$, therefore the result follows.
\end{proof}

\begin{ack} The authors want to thank Professor Noemi Wolanski for providing the proof of Lemma \ref{q1}.
\end{ack}

\appendix
\renewcommand{\theequation}{\Alph{section}.\arabic{equation}}

\section{Properties of $G$ and Orlicz spaces}\label{appAA}
\setcounter{equation}{0} The following results are all include in
\cite{MW1}.
\begin{lema}\label{prop} The function $g$
satisfies the following properties,
\begin{enumerate}
\item[(g1)] $\di \min\{s^{\delta},s^{g_0}\} g(t)\leq g(st)\leq
\max\{s^{\delta},s^{g_0}\} g(t)$
\smallskip
 \item[(g2)] $G$ is convex and $C^2$\item[(g3)]
$\di\frac{t g(t)}{1+g_0}\leq G(t)\leq t g(t) \quad \forall\ t\geq
0.$

\end{enumerate}

\end{lema}
\begin{lema}\label{prop2}
If $\widetilde{G}$ is such that ${\widetilde{G}}'(t)=g^{-1}(t)$
then,
\begin{equation}\label{Gmono1}\tag{$\widetilde{G}1$}
\di\frac{(1+\delta)}{\delta}\min\{s^{1+1/\delta},s^{1+1/g_0}\}\widetilde{G}(t)
\leq \widetilde{G}(st)\leq
\frac{\delta}{1+\delta}\max\{s^{1+1/\delta},s^{1+1/g_0}\}
\widetilde{G}(t)\end{equation}

\end{lema}

We recall that the functional
$$\|u\|_G=\inf\Big\{k>0:\int_{\Omega} G\Big(\frac{|u(x)|}{k}\Big)\,
dx \leq 1\Big\}$$ is a norm in the Orlicz space $L_{G}(\Omega)$
which is the linear hull of the Orlicz class
$$K_G(\Omega)=\Big\{u \mbox{ measurable }:\ \int_{\Omega}
G(|u|)\, dx<\infty\Big\},$$ observe that this set is convex, since
$G$ is also convex (property (g2)). The Orlicz-Sobolev space
$W^{1,G}(\Omega)$ consists of those functions in $L^{G}(\Omega)$
whose distributional derivatives $\nabla u$ also belong to
$L^G(\Omega)$. And we have that
$\|u\|_{W^{1,G}}=\max\{\|u\|_G,\|\nabla u\|_G\}$ is a norm for
this space.

\begin{teo}\label{imb}
 $L^G(\Omega)\hookrightarrow L^{1+\delta}(\Omega)$
continuously.
\end{teo}
\begin{lema}\label{equi}
There exists a constant $C=C(g_0,\delta)$ such that,
$$\|u\|_G \leq C \max\Big\{\Big(\int_{\Omega} G(|u|)\, dx\Big)^{1/(\delta+1)},
\Big(\int_{\Omega} G(|u|)\, dx\Big)^{1/{(g_0+1)}}\Big\}$$
\end{lema}
\begin{lema}\label{poinc} If $u\in W^{1,1}(\Omega)$ with $u=0$ on $\partial\Omega$ and
 $\int_{\Omega} G(|\nabla u|)\,
dx$ is finite, then
$$\int_{\Omega} G\Big(\frac{|u|}{R}\Big)\, dx \leq
\int_{\Omega} G(|\nabla u|)\, dx \quad \mbox{ for
}R=\mbox{diam\,}\Omega.$$
\end{lema}

\section{A result on $\mathcal{L}$-solutions functions with linear growth}\label{appA}
\setcounter{equation}{0}

In this section we will state some properties of
$\mathcal{L}$-subsolutions. From now on, we note $B_r^+=
B_{r}(0)\cap\{x_N>0\}$.

\begin{remark}\label{T} When $|\nabla u|\geq c$, $u$ satisfies  a linear
nondivergence uniformly elliptic equation, $Tu=0$. In our case
\begin{equation}\label{opT}Tv=b_{ij}(\nabla u) D_{i j}v=0\end{equation} where
 $$b_{ij}=\delta_{ij}+\left(\frac{g'(|\nabla u|)|\nabla
 u|}{g(|\nabla u|)}-1\right) \frac{D_i u D_j u}{|\nabla u|^2},$$ and the matrix $b_{ij}(\nabla u)$ is
 $\beta$-elliptic in $\{|\nabla
u|>c\}$, where
 $\beta=\max\{\max\{g_0,1\},\max\{1,1/\delta\}\}$.
\end{remark}
\begin{teo}\label{psub0}
Let $u$ be a Lipschitz function in $\RR^N$ such that
\begin{enumerate}
\item $u\geq 0$ in $\RR^N$, $\mathcal{L} u=0$ in $\{u>0\}$.

\item $\{x_N<0\}\subset \{u>0\}$, $u=0$ in $\{x_N=0\}$.

\item There exists $0<\lambda_0<1$ such that
$\displaystyle\frac{|\{u=0\}\cap B_R(0)|}{|B_R(0)|}>\lambda_0$,
$\forall R>0$.
\end{enumerate}
Then $u=0$ in $\{x_N>0\}$.
\end{teo}
\begin{proof}
See Appendix in \cite{MW2}.
\end{proof}

\section{Blow-up limits}\label{appB}
\setcounter{equation}{0}

Now we  give the definition of blow-up sequence, and we collect
some properties of the limits of these blow-up sequences for
certain classes of functions that are used throughout the paper.

Let $u$ be a function with the following properties,
\begin{enumerate}
\item[(C1)] $u$ is Lipschitz in $\Omega$ with constant $L>0$,
$u\geq 0 \mbox{ in } \Omega$ and $\mathcal{L} u=0 \mbox{ in }
\Omega\cap\{u>0\}$. \item[(C2)] Given $0<\kappa<1$, there exist
two positive constants $C_{\kappa}$ and $r_{\kappa}$ such that for
every ball $B_r(x_0)\subset\Omega$ and $0<r<r_{\kappa}$,
$$
\frac{1}{r}\left(\pint_{B_r(x_0)} u^{\gamma}\, dx
\right)^{1/\gamma}\leq C_{\kappa} \mbox{ implies that } u\equiv 0
\mbox{ in } B_{\kappa r}(x_0).
$$
\item[(C3)] There exist constants $r_0>0$ and
$0<\lambda_1\leq\lambda_2<1$ such that, for every ball
$B_r(x_0)\subset\Omega$ $x_0\mbox{ on }
\partial\{u>0\}$ and $0<r<r_0$
$$
\lambda_1\le
\frac{\left|B_r(x_0)\cap\{u>0\}\right|}{\left|B_r(x_0)\right|} \le
\lambda_2.
$$
\end{enumerate}

\begin{defi}
Let $B_{\rho_k}(x_k)\subset\Omega$ be a sequence of balls with
$\rho_k\to 0$, $x_k\to x_0\in \Omega$ and $u(x_k)=0$. Let
$$
u_k(x):=\frac{1}{\rho_k} u(x_k+\rho_k x).
$$
We call $u_k$ a blow-up sequence with respect to
$B_{\rho_k}(x_k)$.
\end{defi}

Since $u$ is locally Lipschitz continuous, there exists a blow-up
limit $u_0:\R^N\to\R$ such that for a subsequence,
\begin{align*}
& u_k\to u_0 \quad \mbox{in} \quad C^\alpha_{\rm loc}(\R^N)\quad
\mbox{for every}\quad 0<\alpha<1,\\
& \nabla u_k\to\nabla u_0\quad *-\mbox{weakly  in}\quad
L^\infty_{\rm loc}(\R^N),
\end{align*}
and $u_0$ is Lipschitz in $\RR^N$ with constant $L$.
\begin{lema}\label{propblowup}
If $u$ satisfies properties {\rm (C1), (C2)} and {\rm (C3)} then,
\begin{enumerate}
\item $u_0\geq 0$ in $\Omega$ and $\mathcal{L} u_0=0$ in
$\{u_0>0\}$

\medskip

\item $\partial\{u_k>0\}\to \partial\{u_0>0\}$ locally in
Hausdorff distance,

\medskip

\item $\chi_{\{u_k>0\}}\to \chi_{\{u_0>0\}}$ in $L^1_{\rm
loc}(\R^N)$,

\medskip

\item If $K\subset\subset \{u_0=0\}$, then $u_k=0$ in $K$ for big
enough $k$,

\medskip

\item If $K\subset\subset \{u_0>0\}\cup \{u_0=0\}^\circ$, then
$\nabla u_k\rightarrow\nabla u_0$ uniformly in $K$,

\medskip

\item There exists a constant $0<\lambda<1$ such that,
$$
\frac{\left|B_R(y_0)\cap\{u_0=0\}\right|}{\left|B_R(y_0)\right|}\geq
\lambda, \quad \forall R>0, \forall y_0\in \partial\{u_0>0\}
$$

\medskip

\item $\nabla u_k\to\nabla u_0$ a.e in $\R^N$,

\medskip

\item If $x_k\in \partial\{u>0\}$, then $0\in
\partial\{u_0>0\}$
\end{enumerate}
\end{lema}

\begin{proof}
The proof follows as in \cite{FBMW1} and \cite{Le}.
\end{proof}


\def\cprime{$'$}
\providecommand{\bysame}{\leavevmode\hbox
to3em{\hrulefill}\thinspace}
\providecommand{\MR}{\relax\ifhmode\unskip\space\fi MR }
\providecommand{\MRhref}[2]{%
  \href{http://www.ams.org/mathscinet-getitem?mr=#1}{#2}
} \providecommand{\href}[2]{#2}

\end{document}